\documentclass{article}
\usepackage{amsmath,amssymb,bbm}
\catcode`\à=\active \defà{\`a} \catcode`\À=\active \defÀ{\`A}
\catcode`\á=\active \defá{\'a} \catcode`\Á=\active \defÁ{\'A}
\catcode`\ä=\active \defä{\"a} \catcode`\Ä=\active \defÄ{\"A}
\catcode`\â=\active \defâ{\^a} \catcode`\Â=\active \defÂ{\^A}
\catcode`\å=\active \defå{{\aa}} \catcode`\Å=\active \defÅ{{\AA}}
\catcode`\ç=\active \defç{\c{c}} \catcode`\Ç=\active \defÇ{\c{C}}
\catcode`\è=\active \defè{\`e} \catcode`\È=\active \defÈ{\`E}
\catcode`\é=\active \defé{\'e} \catcode`\É=\active \defÉ{\'E}
\catcode`\ë=\active \defë{\"e} \catcode`\Ë=\active \defË{\"E}
\catcode`\ê=\active \defê{\^e} \catcode`\Ê=\active \defÊ{\^E}
\catcode`\ì=\active \defì{\`{\i}} \catcode`\Ì=\active \defÌ{\`{\I}}
\catcode`\í=\active \defí{\'{\i}} \catcode`\Í=\active \defÍ{\'{\I}}
\catcode`\ï=\active \defï{\"{\i}} \catcode`\Ï=\active \defÏ{\"{\I}}
\catcode`\î=\active \defî{\^{\i}} \catcode`\Î=\active \defÎ{\^{\I}}
\catcode`\ò=\active \defò{\`o} \catcode`\Ò=\active \defÒ{\`O}
\catcode`\ó=\active \defó{\'o} \catcode`\Ó=\active \defÓ{\'O}
\catcode`\ö=\active \defö{\"o} \catcode`\Ö=\active \defÖ{\"O}
\catcode`\ô=\active \defô{\^o} \catcode`\Ô=\active \defÔ{\^O}
\catcode`\ù=\active \defù{\`u} \catcode`\Ù=\active \defÙ{\`U}
\catcode`\ú=\active \defú{\'u} \catcode`\Ú=\active \defÚ{\'U}
\catcode`\ü=\active \defü{\"u} \catcode`\Ü=\active \defÜ{\"U}
\catcode`\û=\active \defû{\^u} \catcode`\Û=\active \defÛ{\^U}
\catcode`\ý=\active \defý{\'y} \catcode`\Ý=\active \defÝ{\'Y}
\catcode`\ÿ=\active \defÿ{\"y} \catcode`\˜=\active \def˜{\"Y}
\catcode`\½=\active \def½{!`}
\catcode`\¾=\active \def¾{?`}
\catcode`\ß=\active \defß{{\ss}}
\newcommand{\keywords}[1]{{{\bf Keywords: }#1}}
\newcommand{\AMSclass}[1]{{{\bf A.M.S. subject classification: }#1}}
\newcommand{\tmem}[1]{{\em #1\/}}
\newtheorem{definition}{Definition}
\newcommand{\mathd}{\mathrm{d}}
\newcommand{\tmop}[1]{\operatorname{#1}}
\newtheorem{proposition}{Proposition}
\newcommand{\assign}{:=}
\newtheorem{theorem}{Theorem}
\newtheorem{lemma}{Lemma}
\newtheorem{varnote}{Note}
\newenvironment{note}{\begin{varnote}\em}{\em\end{varnote}}

\begin{document}

\title{Quantizations of the Witt algebra and of simple Lie algebras in
characteristic $p$} \author{Cyril Grunspan\\ \\Dipartimento di Matematica\\Instituto Guido Castelnuovo\\Universita di Roma La Sapienza\\P.le Aldo Moro, 2\\00185 Roma, Italia\\ \\ Email: grunspan[AT]mat.uniroma1.it} \maketitle

\begin{abstract}
  We quantize explicitly the Witt algebra in characteristic $0$ equipped with
  its  Lie bialgebra structures discovered by Taft. Then, we study the
  reduction modulo $p$ of our formulas. This gives $p - 1$ families of
  polynomial noncocommutative deformations of a restricted enveloping algebra
  of a simple Lie algebra in charactristic $p$ (of Cartan type). In
  particular, this yields new families of noncommutative noncocommutative Hopf
  algebras of dimension $p^p$ in char $p$.
  
  {\keywords{Witt algebra, Lie bialgebra, finite dimensional Hopf algebra in
  non-zero characteristic.}}
  
  {\AMSclass{16W30, 17B37, 17B62, 17B66}}
\end{abstract}

{\tmem{``Il est faux de croire que rien n'arrive jamais. En vérité, tout
arrive toujours mais seulement quand on ne le désire plus.''}}

\section{Introduction}

This article comes within the scope of the general program of quantization of
simple Lie bialgebras in non-zero characteristic.

In the following, we use the general method of quantization by Drinfeld's
twist to quantize explicitly the Lie bialgebra structures discovered by Taft
on the Witt algebra in characteristic zero {\cite{T}}. Then, we study the case
of the Witt algebra in characteristic $p$ where $p$ is a prime number and we
show that our formulas lead to quantizations of its restricted enveloping
algebra.

We start by recalling the definitions of the Witt algebra in characteristic
$0$ and of its natural structure of triangular Lie bialgebra of Taft. In
Section $2$, by using the twist discovered by Giaquinto and Zhang {\cite{GZ}},
we quantize explicitly this structure. In Section $3$, we remark that although
the twist used for the quantization involves terms with {\tmem{negative}}
$p$-adic valuation, conjugation by this twist preserves the nonnegative
$p$-adic valuation part of the tensor square of the enveloping algebra of the
Witt algebra. Moreover, the reduction modulo $p$ of these formulas is
compatible with the structure of $p$-Lie algebra on the Witt algebra. We use
this fact to equip the restricted enveloping algebra of the Witt algebra with
non-commutative and non-cocommutative Hopf algebra structures. We do not
expect these Hopf algebras to be triangular since the twist which helps to
define it is not defined on $\mathbbm{F}_p$. We thus get new examples of Hopf
algebra of dimension $p^p$ in characteristic $p$. These Hopf algebras contain
the Radford algebra as a Hopf-subalgebra.

\subsection{Definition of the Witt algebra}

\begin{definition}
  One denotes by $W$ the $\mathbbm{Q}$-Lie algebra given by generators: $L_r,
  r \in \mathbbm{Z}$ and relations: $[ L_r, L_s ] = ( s - r ) L_{r + s}$.
\end{definition}

This Lie algebra is isomorphic to the Lie algebra of vector fields on the
circle and generators $L_r, r \in \mathbbm{Z}$ correspond to operators $x^{r +
1}  \frac{\mathd}{\mathd x}$ in $\tmop{Der} (\mathbbm{Q}[ x, x^{- 1} ] )$. The
following proposition is due to Taft {\cite{T}}.

\begin{proposition}
  There is a triangular Lie bialgebra structure on $W$ given by the $r$-matrix
  $r - r^{21}$ with $r \assign L_0 \otimes L_i$ and $i \in \mathbbm{Z}$.
\end{proposition}

In fact, it has been proved by Ng and Taft that it was the only Lie bialgebra
structure on the ``positive'' part of the Witt algebra {\cite{NgT}}. The
purpose of this article is to quantize this structure and to investigate the
case of the non-zero characteristic.

{\subsection*{Notations.}} In Section $2$, one fixes $i \in \mathbbm{Z}$. In
Section $3$, $p$ is an odd prime number. We set:
\begin{eqnarray*}
  h & \assign & \frac{1}{i} L_0\\
  e & \assign & iL_i
\end{eqnarray*}
For all element $x$ of an unitary $R$-algebra ($R$ a ring) and $a \in R$, we
also set:
\begin{equation}
  x^{( n )} \assign x ( x + 1 ) \ldots ( x + n - 1 ) \label{xap}
\end{equation}

\section{\label{secq}Quantification of Taft's structure}

Let us denote by $( U ( W ), m, \Delta_0, S_0, \varepsilon )$ the natural Hopf
algebra structure on $U ( W )$ i.e.,
\begin{eqnarray*}
  \Delta_0 ( L_k ) & = & L_k \otimes 1 + 1 \otimes L_k\\
  S_0 ( L_k ) & = & - L_k\\
  \varepsilon ( L_k ) & = & 0.
\end{eqnarray*}
In order to quantize the Taft's structure, we recall the twist $F$ of
{\cite{GZ}} defined by:
\[ F \assign \sum_{r = 0}^{\infty}  \frac{1}{r!} h^{( r )} \otimes e^r t^r .
\]
where $t$ denotes a formal variable.

The following proposition has been proved in {\cite{GZ}}.

\begin{proposition}
  We have:
  \begin{eqnarray*}
    ( \Delta_0 \otimes \tmop{Id} ) ( F ) \cdot ( 1 \otimes F ) & = & ( F
    \otimes 1 ) \cdot ( \tmop{Id} \otimes \Delta_0 ) ( F )\\
    ( \varepsilon \otimes \tmop{Id} ) ( F ) & = & ( \tmop{Id} \otimes
    \varepsilon ) ( F ) = 1\\
    F & = & 1 \quad ( t )
  \end{eqnarray*}
  In other words, $F$ is a twist for $( U ( W ), m, \Delta_0, S_0, \varepsilon
  ) i.e., F$ is invertible, $u \assign m \circ ( S \otimes I_d ) ( F )$ is
  invertible and $( U ( W ) \left[ \left[ t \right] \right], m, \Delta_{},
  S_{}, \varepsilon )$ is a Hopf algebra with $\Delta \assign F^{- 1} \Delta_0
  F$ and $S \assign u^{- 1} S_0 u$.
\end{proposition}

We note that $F = 1 + r \cdot t \quad ( t^2 )$. Therefore, the Hopf algebra
constructed with $F$ quantizes the Taft's structure. The following theorem
gives explicitly the quantization.

\begin{theorem}
  \label{to}There exists a structure of non-commutative and non-cocommutative
  Hopf algebra on $U ( W ) \left[ \left[ t \right] \right]$ denoted by $( U (
  W ) \left[ \left[ t \right] \right], m, \Delta_{}, S_{}, \varepsilon )$
  which leaves the product of $U ( W ) \left[ \left[ t \right] \right]$
  undeformed but with a deformed comultiplication defined by:
  \begin{eqnarray}
    \Delta ( L_k ) & \assign & L_k \otimes ( 1 - et )^{\frac{k}{i}} + \sum_{l
    = 0}^{\infty} ( - 1 )^l i^l  \frac{\prod_{j = - 1}^{l - 2} ( k + ji )}{l!}
    h^{( l )} \otimes ( 1 - et )^{- l} L_{k + li} t^l  \label{dellk}\\
    S ( L_k ) & \assign & - ( 1 - et )^{- \frac{k}{i}}  \sum_{l = 0}^{\infty}
    i^l  \frac{\prod_{j = - 1}^{l - 2} ( k + ji )}{l!} L_{k + li} ( h + 1 )^{(
    l )} t^l  \label{slk}\\
    \varepsilon ( L_k ) & \assign & 0  \label{elk}
  \end{eqnarray}
  with $k \in \mathbbm{Z}$.
\end{theorem}

One can simplify formulas of Theorem 4 by introducing the operator $d^{( l )}$
($l \in \mathbbm{N}$) on $U ( W )$ defined by $d^{( l )} \assign \frac{1}{l!}
\tmop{ad} ( e )^l$. Indeed, it is easy to see that
\begin{equation}
  d^{( l )} ( L_k ) = i^l  \frac{\prod_{j = - 1}^{l - 2} ( k + ji )}{l!} L_{k
  + li} \label{dl}
\end{equation}
We have then:
\begin{eqnarray}
  \Delta ( L_k ) & \assign & L_k \otimes ( 1 - et )^{\frac{k}{i}} + \sum_{l =
  0}^{\infty} ( - 1 )^l h^{( l )} \otimes ( 1 - et )^{- l} d^{( l )} ( L_k )
  t^l \\
  S ( L_k ) & \assign & - ( 1 - et )^{- \frac{k}{i}}  \sum_{l = 0}^{\infty}
  d^{( l )} ( L_k ) \cdot ( h + 1 )^{( l )} t^l . 
\end{eqnarray}
We can also give a general formula for the antipode.

\begin{proposition}
  \label{props0}For any homogeneous element $x \in U ( W )$ with respect to
  the graduation given by $|L_k | = k$, one has:
  \begin{equation}
    S ( x ) = ( 1 - et )^{- \frac{|x|}{i}}  \sum_{n = 0}^{\infty} d^{( n )} (
    S_0 ( x ) ) \cdot ( h + 1 )^{( n )} t^n .
  \end{equation}
\end{proposition}

In fact, by applying Lemma \ref{grazie} below with $a = i$ and $k$ replaced by
$k - i$, we remark that the structure coefficients of $\Delta$ and $S$ in
(\ref{dellk}) and (\ref{slk}) belong to $\mathbbm{Z}$.

\begin{lemma}
  \label{grazie}For all integers $a, k$ and $l$, $a^l  \frac{\prod_{j = 0}^{l
  - 1} ( k + ja )}{l!}$ is an integer.
\end{lemma}

This fact allows us to consider the reduction modulo $p$ of our formulas.

\section{Quantizations of restricted enveloping algebras}

We show that formulas above can be used to define an Hopf algebra structure on
the {\tmem{restricted}} enveloping algebra of the Witt algebra in
characteristic $p$ {\cite{J}}.

\begin{definition}
  One denotes by $\mathfrak{D}$ the $\mathbbm{F}_p$-Lie algebra given by
  generators: $D_k, k \in \mathbbm{Z}$ and relations:
  \begin{eqnarray*}
    [ D_k, D_l ] & = & ( l - k ) D_{k + l}\\
    D_{k + p} & = & D_k
  \end{eqnarray*}
  for $k, l \in \mathbbm{Z}$. The algebra $\mathfrak{D}$ is called the Witt
  algebra in characteristic $p$.
\end{definition}

The algebra $\mathfrak{D}$ is isomorphic to $\tmop{Der} (\mathbbm{F}_p [ X ] /
( X^p - 1 ) )$ and $D_k$ corresponds to the operator $X^{k + 1} 
\frac{\mathd}{\mathd X}$.

It is well known that $\mathfrak{D}$ is a simple restricted Lie algebra
{\cite{J}}. Its structure of $p$-Lie algebra is given by $D_0^{( 0 )} = D_0$
and $D_k^{( 0 )} = 0$ for $k$ non divisible by$p$. Its restricted enveloping
algebra $\text{$U_c (\mathfrak{D})$}$ is isomorphic to $U (\mathfrak{D}) / I$
where $I$ is the ideal of $U (\mathfrak{D})$ generated by $D_r^p - D_r$ with
$p$ divides $r$ and $D_k^p$ with $k$ non divisible by $p$. A basis of this
algebra is given by monomials $\prod_{k = 0}^{p - 1} D_k^{\alpha_k}$ with
$\alpha_k \in \{ 0, \ldots, p - 1 \}$ {\cite{J}}. So, $\dim_{\mathbbm{F}_p}
U_c (\mathfrak{D}) = p^p$.

\begin{note}
  The Witt algebra in characteristic $p$ is also sometimes defined as the Lie
  algebra $\mathfrak{w} \mathfrak{i} \mathfrak{t} \mathfrak{t}$ given by
  generators: $e_k, k \in \{ - 1, \ldots, p - 2 \}$ and relations:
  \begin{equation}
    [ e_k, e_l ] = \left\{\begin{array}{cc}
      ( l - k ) e_{k + l} & \text{if } k + l \leqslant p - 2 ;\\
      0 & \text{otherwise} .
    \end{array}\right. \label{ekel}
  \end{equation}
  The Lie algebra $\mathfrak{w} \mathfrak{i} \mathfrak{t} \mathfrak{t}$ is
  isomorphic to $\tmop{Der} (\mathbbm{F}_p [ Y ] / ( Y^p ) )$ and $e_k$ stands
  for the operator $Y^{k + 1}  \frac{\mathd}{\mathd Y}$. This Lie algebra is
  also isomorphic to $\mathfrak{D}$ as it can be seen by setting $X = Y + 1$.
  The corresponding Lie algebra morphism maps $e_k$ to $\sum_{l = - 1}^k ( - 1
  )^l  \left(\begin{array}{c}
    k + 1\\
    l + 1
  \end{array}\right) D_l$.
\end{note}

On $\mathbbm{F}_p$, the twist $F$ does not make sense because of the
coefficient $\frac{1}{n!}$ in the definition of $F$. However, as we have seen
at the end of Section \ref{secq}, formulas of Theorem \ref{to} are defined on
$\mathbbm{Z}$. Moreover, Lemma \ref{compp} below shows that these coefficients
are compatible with the reduction modulo $p$.

\begin{lemma}
  \label{compp}Let $\alpha, \kappa, l$ be integers. Then, the residue class of
  $\alpha^l  \frac{\prod_{j = 0}^{l - 1} ( \kappa + j \alpha )}{l!}$ modulo
  $p$ depends only on $l$ and on the residue class of $\alpha$ and $\kappa$.
\end{lemma}

By definition, for $l \in \mathbbm{N}$ and $a, k \in \mathbbm{F}_p$, we will
denote by $N ( a, k, l )$ the common residue modulo $p$ of all integers of the
form $\alpha^l  \frac{\prod_{j = 0}^{l - 1} ( \kappa + j \alpha )}{l!}$ with
$\tmop{cl} ( \alpha ) = a, \tmop{cl} ( \kappa ) = k$.

Therefore, Theorem \ref{tq} below make sense in $\mathbbm{F}_p$.

\begin{theorem}
  \label{tq}Let $i$ be an element of $\mathbbm{F}_p - \{ 0 \}$. Then, the Hopf algebra $( U_c (\mathfrak{D}) \left[ t \right], m, \Delta_{},
  S_{}, \varepsilon )$ is a
  polynomial deformation  of the restricted enveloping algebra of $\mathfrak{D}$.
  The algebra structure is undeformed and the coalgebra structure is given by:
  \begin{eqnarray}
    \Delta ( D_k ) & \assign & D_k \otimes ( 1 - et )^{\frac{k}{i}} + \sum_{l
    = 0}^{p - 1} ( - 1 )^l N ( i, k - i, l ) h^{( l )} \otimes ( 1 - et )^{-
    l} D_{k + li} t^l  \label{deldk}\\
    S ( D_k ) & \assign & - ( 1 - et )^{- \frac{k}{i}}  \sum_{l = 0}^{p - 1} N
    ( i, k - i, l ) D_{k + li} ( h + 1 )^{( l )} t^l  \label{sdk}\\
    \varepsilon ( D_k ) & \assign & 0  \label{edk}
  \end{eqnarray}
  with, $k \in \mathbbm{F}_p, h \assign \frac{1}{i} D_0$ and $e \assign iD_i$.
\end{theorem}

The above sums are finite, which explains why they define a polynomial
deformation. Thus, we can specialize $t$ to any element of $\mathbbm{F}_p$.
Since $i$ is an arbitrary element of $\mathbbm{F}_p - \{ 0 \}$, this gives $p
- 1$ new families of non-commutative and non-cocommutative Hopf algebra of
dimension $p^p$ in characteristic $p$.

\begin{note}
  If we set $\alpha \assign ( 1 - et )^{- 1}$, then we have: $\left[ h, \alpha
  \right] = \alpha^2 - \alpha ; h^p = h ; \alpha^p = 1 ; \Delta ( h ) = h
  \otimes \alpha + 1 \otimes h, \alpha$ is ``group-like'', $S ( h ) = h
  \alpha^{- 1}$ and $\varepsilon ( h ) = 0$.
  
  So, the sub-algebra generated by $h$ and $e$ is a sub-Hopf algebra of $U_c
  (\mathfrak{D})$ isomorphic to the Radford algebra .
\end{note}

{\section*{Acknowledgement}}

I am very grateful to B. Enriquez for help and support during the preparation
of this article. I am also grateful to P. Bauman, J. Bichon, F. Gavarini and
E. Taft for comments and S-H Ng who kindly emailed me his thesis. I would also
like to thank C. Procesi for the hospitality of the Castelnuovo Institute and
C. Kassel for the hospitality of the university of Strasbourg.

\bibliographystyle{plain}

\end{document}